\documentclass[11pt]{article}
\usepackage[T1]{fontenc}
\usepackage{lmodern}
\usepackage[margin=1in]{geometry}
\usepackage{amsmath,amssymb,amsthm,mathtools,mathrsfs}
\usepackage{microtype,enumitem,xcolor,xurl}
\usepackage{hyperref}
\hypersetup{colorlinks=true,linkcolor=teal!65!black,citecolor=teal!65!black,
urlcolor=teal!65!black,pdfauthor={Yi Wang and Shuhan Yang},
pdftitle={Smooth Exterior Potentials for a Gradient-Weighted Hessian Equation}}
\numberwithin{equation}{section}
\allowdisplaybreaks[2]
\newtheorem{theorem}{Theorem}[section]
\newtheorem{lemma}[theorem]{Lemma}
\newtheorem{proposition}[theorem]{Proposition}
\newcommand{\RR}{\mathbb R}
\newcommand{\Sphere}{\mathbb S}
\newcommand{\Div}{\operatorname{div}}

\newcommand{\ip}[2]{\langle #1,#2\rangle}
\newcommand{\LL}{\mathcal L_{k,p}}
\title{Smooth Exterior Potentials for a Gradient-Weighted Hessian Equation}
\author{Yi Wang\thanks{Department of Mathematics, Johns Hopkins University,
3400 N Charles St, Baltimore, MD 21218, USA.
\href{mailto:ywang261@jhu.edu}{\mbox{\textcolor{black}{\texttt{ywang261@jhu.edu}}}}}
\and Shuhan Yang\thanks{Department of Mathematics, Johns Hopkins University,
3400 N Charles St, Baltimore, MD 21218, USA.
\href{mailto:syang180@jhu.edu}{\mbox{\textcolor{black}{\texttt{syang180@jhu.edu}}}}}}
\date{November 22, 2026}
\begin{document}
\maketitle
\begin{abstract}
We construct a positive solution of an exterior Dirichlet problem for a gradient-weighted
$k$-Hessian operator on each smooth strictly star-shaped domain whose boundary
curvatures lie in $\Gamma_k$. For $p>1$ sufficiently close to one, the solution
is smooth up to the boundary, has no critical points, and is unique among
continuous potentials on the natural viscosity branch. Its levels remain
strictly star-shaped and strictly $k$-convex. We determine its radial
asymptotics through two derivatives and the uniform limit of $(p-1)|Du_p|$
on the boundary. The construction perturbs an inverse curvature-quotient
flow. In logarithmic radial coordinates the exact equation is elliptic for
$p>1$ and its limit is parabolic. A weighted inverse estimate on the
half-cylinder controls the singular perturbation.
\end{abstract}

\section{Introduction and main result}

The exterior $p$-capacitary potential connects the geometry of a domain
to an equation on its complement. For convex domains its levels satisfy
geometric monotonicity formulas \cite{FMP2019}. Its far-field behavior
on asymptotically conical manifolds was studied in \cite{BFMAsymptotics2024}.
In the proof of Agostiniani,
Fogagnolo and Mazzieri \cite{AFM}, estimates for the flux of this
potential yield a Minkowski inequality as $p$ decreases to one.
To measure higher curvatures, one can replace the scalar flux
$|Du|^{p-2}Du$ by a flux involving the $(k-1)$st Newton tensor of
the Hessian. This raises an existence question. A boundary that is
strictly $k$-convex may have negative principal curvatures, and the
operator depends on the direction of the potential gradient. We
construct a potential for which this direction is defined everywhere
and the required curvature cone persists through all exterior levels.

Write $E=\RR^n\setminus\overline\Omega$. For a symmetric matrix $A$,
let $\sigma_j(A)$ be the $j$th elementary symmetric polynomial of its
eigenvalues, with $\sigma_0=1$. Set $T_{-1}=0$, $T_0=I$ and
$T_j(A)=\sigma_j(A)I-AT_{j-1}(A)$. We use the cone
$\Gamma_k=\{\lambda:\sigma_j(\lambda)>0,\ 1\le j\le k\}$ in the
dimension of the matrix under discussion. The equation and its data are
\begin{equation}\label{eq:exterior}
 \LL u:=\Div\bigl(|Du|^{p-k-1}T_{k-1}(D^2u)Du\bigr)=0
 \quad\text{in }E,\qquad u|_{\partial\Omega}=1,\qquad
 \lim_{|x|\to\infty}u(x)=0 .
\end{equation}
For $k=1$, \eqref{eq:exterior} is the $p$-Laplace equation. At $p=k+1$,
the weight is constant and $z=-u$ solves $\sigma_k(D^2z)=0$ in $E$,
with $z=-1$ on $\partial\Omega$ and $z\to0$ at infinity.
For Hessian equations on bounded domains, Tso symmetrized convex
sublevel sets while preserving their $(k-1)$st quermassintegral
\cite{TSO89}. For $1\le k<n/2$
and a smooth star-shaped domain with strictly $(k-1)$-convex boundary,
Ling Xiao constructed an admissible $C^{1,1}$ exterior solution with regular
levels and used its radial asymptotics to prove a weighted Minkowski
inequality \cite[Theorem~1.2 and Corollary~1.1]{Xiao2022}. For convex
domains in this class, Ma and Zhang proved uniqueness among $k$-convex
solutions and obtained two-sided gradient bounds and a Hessian bound at
infinity \cite[Theorem~1.1]{MaZhang2024}. Gao, Ma and Zhang studied the
real homogeneous equation on bounded punctured domains \cite{GaoMaZhang2023}
and the homogeneous complex $k$-Hessian equation on exterior domains
\cite{Gao2023Exterior}. Li and Xiao treated the exterior equation
$\sigma_k(D^2z)=1$ with prescribed quadratic asymptotics
\cite[Theorem~1.4]{LiXiao2026}.

We consider $p>1$ close to one. The Hessian is then stretched in the
gradient direction. At a nonzero vector $\xi$, put
\begin{equation}\label{eq:weight}
 P(\xi)=I+\left(\sqrt{(p-1)/k}-1\right)
               \frac{\xi\otimes\xi}{|\xi|^2},\qquad
 M[u]=-P(Du)D^2uP(Du).
\end{equation}
The positive-potential convention uses $M[u]\in\partial\Gamma_k$.
For a continuous potential this means that each smooth test $\phi$
touching from above at an interior point with $D\phi\ne0$ has
$M[\phi]\notin\Gamma_k$, and each such test from below has
$M[\phi]\in\overline\Gamma_k$. These are the nonzero-gradient
branch tests used in the uniqueness assertion.

We orient each level by $-Du/|Du|$ and take the outer shape operator
$h_\delta=D\nu$ on tangent vectors, so that spheres have positive
curvatures. Write $H_j=\sigma_j(h_\delta)/\binom{n-1}{j}$.
After a translation, strict star-shapedness means that
$\partial\Omega=\{r_\Omega(\omega)\omega:\omega\in\Sphere^{n-1}\}$
for a positive smooth function $r_\Omega$. Strict $k$-convexity means
that the boundary curvature vector belongs to $\Gamma_k$ in dimension $n-1$.

\begin{theorem}\label{thm:starshaped}
Let $n\ge3$, $1\le k\le n-2$, and let $\Omega$ be a bounded smooth
strictly star-shaped, strictly $k$-convex domain. There exists
$0<\varepsilon_0(\Omega,n,k)<n-k-1$ such that whenever
$1<p<1+\varepsilon_0$, \eqref{eq:exterior} has a solution $u_p$ with
the following properties.
\begin{enumerate}[label=\textup{(P\arabic*)},leftmargin=*]
\item $u_p$ is smooth up to the finite boundary, with $0<u_p<1$ in $E$.
\item The gradient is nonzero throughout $\overline E$, the equation
holds classically, and $M[u_p]\in\partial\Gamma_k$ in dimension $n$.
\item Every level in $0<u_p\le1$ is a smooth strictly star-shaped,
strictly $k$-convex hypersurface. In particular $\sigma_{k-1}(h_\delta)>0$.
\item For $z=-u_p$, the outward Newton field
$|Dz|^{p-k-1}T_{k-1}(D^2z)Dz$ is smooth and divergence-free in $E$,
also in the sense of distributions.
\item For some $\mathsf a_p>0$ and $\alpha_p=(n-k-p+1)/(p-1)$,
uniformly in direction,
\[
 D^j\bigl(u_p(x)-\mathsf a_p|x|^{-\alpha_p}\bigr)
 =o(|x|^{-\alpha_p-j}),\qquad j=0,1,2.
\]
\end{enumerate}
This solution is unique among positive continuous exterior potentials
satisfying the data in \eqref{eq:exterior} and the branch tests above.
If $B_R(x_0)\subset\Omega$, then $(R/|x-x_0|)^{\alpha_p}\le u_p(x)$ on $\overline E$.
On the fixed boundary,
\begin{equation}\label{eq:boundarylimit}
 (p-1)|Du_p|\ \longrightarrow\ k\frac{\sigma_k(h_\delta)}{\sigma_{k-1}(h_\delta)}
   =(n-k)\frac{H_k}{H_{k-1}}\qquad\text{uniformly as }p\downarrow1.
\end{equation}
\end{theorem}

The boundary limit is the link with higher-order Minkowski inequalities.
Since $H_k/H_{k-1}$ has a positive minimum on $\partial\Omega$,
\eqref{eq:boundarylimit} implies $(|Du_p|/\alpha_p)^{p-1}\to1$
uniformly there. Thus the boundary integral of
$(|Du_p|/\alpha_p)^{p-1}H_{k-1}$ tends to the curvature integral
$\int_{\partial\Omega} H_{k-1}\,\mathrm{d}S_\delta$ in the consecutive
Minkowski inequality of Guan and Li \cite{GuanLi2009}.

To construct the potential, set $\varepsilon=p-1$ and parametrize its
levels by $s=-\varepsilon\log u$. For $k=1$ this gives
$\Div(|Ds|^{p-2}Ds)=|Ds|^p$. Its formal limit where $Ds\ne0$ is
$\Div(Ds/|Ds|)=|Ds|$, the level-set equation for inverse mean curvature
flow. Huisken and Ilmanen developed the weak theory of this flow
\cite{HuiskenIlmanen2001}, and Moser constructed weak solutions on
Euclidean exterior domains using $p$-harmonic functions \cite{Moser2007}.
For higher $k$, the radial equation at $\varepsilon=0$ gives the outward
normal speed $\sigma_{k-1}/(k\sigma_k)$. Gerhardt and Urbas studied smooth
expanding flows driven by inverse curvature functions
\cite{Gerhardt1990,Urbas1990}. For this particular quotient, the theorem of
Ding and Li \cite{DingLi2022} supplies the global reference flow and its
exponential convergence after rescaling (Proposition~\ref{star:reference}).
VanBlargan and Wang used the rescaled quotient flow to study the stability
of the quermassintegral inequality near spheres \cite{VanBlarganWangFlow2024}.

The article is organized as follows. In Section~\ref{sec:preliminaries}
we prove the power comparison principle, derive the radial equation,
and linearize it along the reference flow. Section~\ref{smsp:linearinverse}
establishes a weighted inverse estimate on the half-cylinder.
In Section~\ref{sec:perturbation} we use this estimate to solve the radial
equation for $p>1$ close to one. In Section~\ref{sec:main-proof} we
reconstruct the exterior potential and prove Theorem~\ref{thm:starshaped}.

\section{Preliminaries}\label{sec:preliminaries}

\subsection{The weighted Hessian and power comparison}

In this subsection we relate the divergence of the Newton flux to the
weighted Hessian and use power comparison to establish uniqueness on the
natural branch and the lower bound supplied by the potential of an inscribed ball.

\begin{samepage}
\begin{lemma}[Newton flux and the ball potential]\label{sec:weightedhessian}
Let $p>1$. For a $C^3$ function $u$ with $Du\ne0$, set $z=-u$ and
$F_z=|Dz|^{p-k-1}T_{k-1}(D^2z)Dz$. Then
\begin{equation}\label{eq:newtonidentity}
 (-1)^k\LL u=\Div F_z=k|Du|^{p-k-1}\sigma_k(M[u])
 =k\sigma_k\!\left(D\!\left(-|Du|^{(p-k-1)/k}Du\right)\right).
\end{equation}
If $p-1<n-k$, then on the exterior of $B_R(x_0)$ the function
$w(x)=(R/|x-x_0|)^{\alpha_p}$ satisfies \eqref{eq:exterior} and
$M[w]\in\partial\Gamma_k$.
\end{lemma}
\end{samepage}

\begin{proof}
Write $A=-D^2u$ and $\nu=-Du/|Du|$. The Newton tensor $T_{k-1}(D^2z)$
is divergence-free: in its alternating-index formula each third derivative
is symmetric in two indices over which the formula is antisymmetric.
Since $T_{k-1}(-D^2u)=(-1)^{k-1}T_{k-1}(D^2u)$ and
$\operatorname{tr}(AT_{k-1}(A))=k\sigma_k(A)$, differentiating the flux gives
\[
 \Div F_z=(-1)^k\LL u
 =|Du|^{p-k-1}\left[k\sigma_k(A)+(p-k-1)
            \ip{AT_{k-1}(A)\nu}{\nu}\right].
\]
In a frame ending with $\nu$, the last contraction equals
$\sigma_k(A)-(T_k(A))_{\nu\nu}$, the sum of the principal $k$-minors
involving the normal index. The congruence $A\mapsto PAP$, with $P=P(Du)$, multiplies precisely
these minors by $(p-1)/k$, giving $\Div F_z=k|Du|^{p-k-1}\sigma_k(PAP)$.
The Jacobian in \eqref{eq:newtonidentity} equals
$|Du|^{(p-k-1)/k}P^2A$. Since $P^2A$ is similar to $PAP$, homogeneity
gives the last equality.

For $d=|x-x_0|$, the tangential eigenvalues of $M[w]$ equal
$\alpha_p w/d^2$, and its radial eigenvalue equals
$-(n-k)\alpha_p w/(kd^2)$. Hence for $1\le j\le k$,
\[
 \sigma_j(M[w])=
 \left(\frac{\alpha_p w}{d^2}\right)^j
 \binom{n-1}{j-1}\left(\frac{n-j}{j}-\frac{n-k}{k}\right).
\]
This is positive for $j<k$ and zero for $j=k$. The function $w$ has
nonzero gradient and the required boundary values, so \eqref{eq:newtonidentity}
gives the equation. At $d=R$ one has $(p-1)|Dw|=(n-k-p+1)/R$,
consistent with \eqref{eq:boundarylimit}.
\end{proof}

Trudinger and Wang \cite[\S4, (4.3)]{TW} use the Jacobian
$D(|Dz|^{r-2}Dz)$ for $r\ge2$ in their estimates for Hessian measures.
In \eqref{eq:newtonidentity} the corresponding exponent is $r=1+(p-1)/k$.

\begin{proposition}[Power comparison]\label{star:boundary}
Let $1<p<n-k$, and suppose $u$ is a positive solution of
\eqref{eq:exterior}, smooth up to the finite boundary, with
$0<u<1$ in $E$ and $Du\ne0$ on $\overline E$. Assume that
$M[u]\in\partial\Gamma_k$ in $E$ and that the levels of $u$ are
strictly $k$-convex with respect to $-Du/|Du|$. Then $u$ is the
unique positive continuous potential with the same boundary data
and limit at infinity satisfying the nonzero-gradient branch tests.
Every ball $B_R(x_0)\subset\Omega$ satisfies
\[
 \left(\frac{R}{|x-x_0|}\right)^{\alpha_p}\le u(x)
 \qquad (x\in\overline E).
\]
At a point where this ball is tangent to $\partial\Omega$,
$|Du|/\alpha_p\le 1/R$.
\end{proposition}

\begin{proof}
An upper test for $u$ has the same gradient and a larger Hessian
at contact, hence a smaller weighted matrix than $M[u]$. Its matrix
cannot belong to $\Gamma_k$, which is preserved by the addition of a
positive semidefinite matrix. Similarly the matrix of a lower test
belongs to $\overline\Gamma_k$. Thus $u$ satisfies both branch tests.

For $\gamma>0$, $D(u^\gamma)$ is a positive multiple of $Du$, so the
two gradients define the same $P$. The chain rule gives
\[
 M[u^\gamma]=\gamma u^{\gamma-1}M[u]
  +\gamma(1-\gamma)u^{\gamma-2}(PDu)\otimes(PDu).
\]
The perturbation changes only the normal diagonal entry in a frame
adapted to the level. The tangential block of $M[u]$ is $|Du|h_\delta$,
and $|PDu|^2=(p-1)|Du|^2/k$. Expansion of the principal minors
therefore gives, for $1\le j\le k$,
\[
 \begin{aligned}
 \sigma_j(M[u^\gamma])
 ={}&(\gamma u^{\gamma-1})^j\sigma_j(M[u])\\
 &+\gamma(1-\gamma)u^{\gamma-2}\frac{p-1}{k}|Du|^2
      (\gamma u^{\gamma-1}|Du|)^{j-1}\sigma_{j-1}(h_\delta).
 \end{aligned}
\]
Strict $k$-convexity makes each $\sigma_{j-1}(h_\delta)$ positive.
Since $M[u]\in\partial\Gamma_k$, its lower symmetric polynomials
are nonnegative and $\sigma_k(M[u])=0$. Consequently $M[u^\gamma]$
belongs to $\Gamma_k$ for $0<\gamma<1$, whereas
$\sigma_k(M[u^\gamma])<0$ for $\gamma>1$.

Let $\widetilde u$ be a potential satisfying the stated tests and
data. A positive maximum of $\widetilde u-u^\gamma$ for $0<\gamma<1$
is attained in $E$, since the difference vanishes on the boundary and
at infinity. Translating $u^\gamma$ up by this maximum gives an upper
test with nonzero gradient and matrix in $\Gamma_k$. This contradicts
the upper test condition. For $\gamma>1$, a positive maximum of
$u^\gamma-\widetilde u$ gives, after translation down, a lower test
with $\sigma_k<0$, contradicting the lower test condition. Taking
$\gamma\to1$ from the two sides gives $\widetilde u=u$.

By Lemma~\ref{sec:weightedhessian}, $w=(R/|x-x_0|)^{\alpha_p}$ solves the ball equation, has
positive level curvatures, satisfies $w\le1$ on $\partial\Omega$,
and tends to zero at infinity. For $\gamma>1$ the same calculation
gives $\sigma_k(M[w^\gamma])<0$. A positive maximum of $w^\gamma-u$
would thus furnish a forbidden lower test for $u$. Letting
$\gamma\downarrow1$ gives $w\le u$. At boundary tangency $w=u=1$,
so differentiation along the common outward normal gives
$-|Du|=\partial_\nu u\ge\partial_\nu w=-\alpha_p/R$.
\end{proof}

\subsection{The reference flow and the radial equation}

In this subsection we describe the limiting inverse curvature flow,
derive the exact equation for the radial graphs of the level sets,
and compute its linearization along the reference flow.

Assume the domain hypotheses of Theorem~\ref{thm:starshaped}. On the
normalized radial graph $r_0(\omega,s)\omega$, denote the outward normal
and shape operator by $\nu_0$ and $h_0$.

\begin{proposition}[Reference flow]\label{star:reference}
The flow from $\partial\Omega$ with normal speed
$\sigma_{k-1}/(k\sigma_k)$ exists for all $s\ge0$. After dilation by
$e^{-s/(n-k)}$, its radial function $r_0$ has initial value $r_\Omega$ and satisfies
$$
(r_0)_s+\frac{r_0}{n-k}
=\frac{1}{k(\nu_0\cdot\omega)}
\frac{\sigma_{k-1}(h_0)}{\sigma_k(h_0)}.
$$
The functions $r_0$, $\nu_0\cdot\omega$ and $(r_0)_s+r_0/(n-k)$
have uniform positive lower bounds. The curvature vectors of $h_0$
remain in a compact subset of $\Gamma_k$. Moreover, $v_0=\log r_0$
converges to a constant $v_{0,\infty}$ exponentially, with all
derivatives of any prescribed finite order in $(\omega,s)$.
\end{proposition}

\begin{proof}
Apply the inverse curvature flow theorem of Ding and Li
\cite[Assumption~1.1 and Theorem~1.2]{DingLi2022}.
Take their curvature function to be $k\sigma_k/\sigma_{k-1}$ and
their parameters to be $\alpha=\delta=0$, $\beta=1$. The quotient is
positive, homogeneous of degree one, concave and strictly increasing
on $\Gamma_k$ by \cite[(1.5)]{ITW}. It extends continuously by zero
on the boundary of the cone, since $H_k/H_{k-1}\le H_k^{1/k}$ by the
Maclaurin inequalities for $k\ge2$, and $H_1=0$ on $\partial\Gamma_1$.
The initial surface is an admissible radial graph, and the curvature
function equals $n-k$ on the unit sphere. The cited theorem gives
a global smooth flow and exponential smooth convergence of its
normalized graphs to a sphere of positive radius centered at the origin.

Taking the normal component of the velocity of
$e^{s/(n-k)}r_0\omega$ gives the displayed radial equation by
homogeneity of the speed. On each bounded time interval the lower
bounds and the compact curvature range follow from strict
star-shapedness, admissibility and compactness. The spherical limit
gives the same bounds on the tail. The radial equation is autonomous
and its constant solutions are stationary, so successive differentiation
in $s$ gives exponential decay for time derivatives as well as angular
derivatives.
\end{proof}

Fix a common decay rate $\rho>0$ for the finitely many derivatives
needed below. Set $\varepsilon=p-1$ and suppose $0<\varepsilon<n-k$.
For a smooth positive radial function $r(\omega,s)$ write
$$
x(\omega,s)=e^{s/(n-k-\varepsilon)}r(\omega,s)\omega,
\qquad u(x(\omega,s))=e^{-s/\varepsilon}.
$$
Let $\nu$, $h$, and $\nabla^\top$ denote the outward normal, shape operator,
and tangential gradient of the normalized surface $r\omega$ in its
Euclidean induced metric. In an orthonormal frame for the unit sphere,
$$
\nu=\frac{r\omega-\nabla_{\mathbb S}r}
{\sqrt{r^2+|\nabla_{\mathbb S}r|^2}},
\qquad
(g_r)_{ij}=r^2\delta_{ij}+r_i r_j.
$$
Thus $\nu\cdot\omega>0$. Define the outward radial speed by
$$
V=r_s+\frac{r}{n-k-\varepsilon}
$$
and require $V>0$. The physical arrival coordinate satisfies
$$
|Ds|=\frac{e^{-s/(n-k-\varepsilon)}}{(\nu\cdot\omega)V}.
$$

\begin{lemma}[Exact radial equation]\label{star:equation}
Suppose $r>0$, $V>0$ and $\sigma_{k-1}(h)>0$. In the coordinates
above, $\LL u=0$ is equivalent to $\mathscr F_\varepsilon(r)=0$, where
\begin{equation}\label{eq:radialoperator}
\begin{aligned}
\mathscr F_\varepsilon(r)={}&\varepsilon r_{ss}
+\left(1+\frac{2\varepsilon}{n-k-\varepsilon}\right)r_s
+\left(\frac1{n-k-\varepsilon}
+\frac{\varepsilon}{(n-k-\varepsilon)^2}\right)r\\
&-k(\nu\cdot\omega)V^2\frac{\sigma_k(h)}{\sigma_{k-1}(h)}
+\varepsilon V\partial_s\log(\nu\cdot\omega)\\
&+\varepsilon(\nu\cdot\omega)V^3
\bigl(\omega-(\nu\cdot\omega)\nu\bigr)\cdot
\nabla^\top\!\left(\frac1{(\nu\cdot\omega)V}\right)\\
&+\varepsilon\frac{(\nu\cdot\omega)^3V^4}{\sigma_{k-1}(h)}
T_{k-2}(h)\left(
\nabla^\top\!\left(\frac1{(\nu\cdot\omega)V}\right),
\nabla^\top\!\left(\frac1{(\nu\cdot\omega)V}\right)
\right).
\end{aligned}
\end{equation}
The coefficient of $r_{ss}$ is $\varepsilon$, and every term involving
an angular-time second derivative has a factor of $\varepsilon$.
\end{lemma}

\begin{proof}
Since $Du=-(u/\varepsilon)Ds$ and $V>0$, the normal in
\eqref{eq:weight} is $Ds/|Ds|=\nu$. With $u=e^{-s/\varepsilon}$, the tangential, mixed and normal blocks of
$-PD^2uP$, divided by the positive factor $u/\varepsilon$, are
\begin{equation}\label{eq:physicalblock}
\begin{pmatrix}
|Ds|h_\delta&\sqrt{\varepsilon/k}\,\nabla^{\top,\delta}|Ds|\\
\sqrt{\varepsilon/k}\,(\nabla^{\top,\delta}|Ds|)^{\mathsf T}&
\bigl(\varepsilon D_\nu|Ds|-|Ds|^2\bigr)/k
\end{pmatrix}.
\end{equation}
Here $h_\delta$ and $\nabla^{\top,\delta}$ are on the physical level.
By \eqref{eq:newtonidentity}, $\LL u=0$ precisely when the sum of the
principal $k$-minors of \eqref{eq:physicalblock} vanishes.
The tangential minors sum to $|Ds|^k\sigma_k(h_\delta)$, the
coefficient of the normal entry is $|Ds|^{k-1}\sigma_{k-1}(h_\delta)$,
and the mixed entries contribute
$-\varepsilon |Ds|^{k-2}
T_{k-2}(h_\delta)(\nabla^{\top,\delta}|Ds|,
\nabla^{\top,\delta}|Ds|)/k$.
Setting their sum to zero and dividing by $|Ds|^{k-2}/k>0$
therefore gives exactly
$$
\begin{aligned}
|Ds|^3\sigma_{k-1}(h_\delta)-k|Ds|^2\sigma_k(h_\delta)
=\varepsilon\bigl[&|Ds|\sigma_{k-1}(h_\delta)D_\nu|Ds|\\
&-T_{k-2}(h_\delta)
(\nabla^{\top,\delta}|Ds|,\nabla^{\top,\delta}|Ds|)\bigr].
\end{aligned}
$$
For $k=1$ the mixed term vanishes since $T_{-1}=0$.

The parameter curves have velocity
$x_s=e^{s/(n-k-\varepsilon)}V\omega$. Their tangential part must be
retained when differentiating $|Ds|$. The chain rule gives
$$
\begin{aligned}
h_\delta&=e^{-s/(n-k-\varepsilon)}h,\\
\nabla^{\top,\delta}|Ds|
&=e^{-2s/(n-k-\varepsilon)}
  \nabla^\top\!\left(\frac1{(\nu\cdot\omega)V}\right),\\
D_\nu|Ds|
&=\frac{e^{-2s/(n-k-\varepsilon)}}{(\nu\cdot\omega)V}
\biggl[\partial_s\!\left(\frac1{(\nu\cdot\omega)V}\right)
       -\frac1{(n-k-\varepsilon)(\nu\cdot\omega)V}\\
&\hspace{42mm}
-V\bigl(\omega-(\nu\cdot\omega)\nu\bigr)\cdot
\nabla^\top\!\left(\frac1{(\nu\cdot\omega)V}\right)\biggr].
\end{aligned}
$$
Substituting and canceling the common positive factors yields
$$
\begin{aligned}
\frac1{(\nu\cdot\omega)V}
-k\frac{\sigma_k(h)}{\sigma_{k-1}(h)}
={}&\varepsilon\biggl[
\partial_s\!\left(\frac1{(\nu\cdot\omega)V}\right)
-\frac1{(n-k-\varepsilon)(\nu\cdot\omega)V}\\
&-V\bigl(\omega-(\nu\cdot\omega)\nu\bigr)\cdot
\nabla^\top\!\left(\frac1{(\nu\cdot\omega)V}\right)\\
&-\frac{(\nu\cdot\omega)^2V^2}{\sigma_{k-1}(h)}
T_{k-2}(h)\left(
\nabla^\top\!\left(\frac1{(\nu\cdot\omega)V}\right),
\nabla^\top\!\left(\frac1{(\nu\cdot\omega)V}\right)
\right)\biggr].
\end{aligned}
$$
Move the right-hand side to the left and multiply by
$(\nu\cdot\omega)V^2>0$. The identity
$$
\partial_s\log\!\left(\frac1{(\nu\cdot\omega)V}\right)
=-\frac{V_s}{V}-\partial_s\log(\nu\cdot\omega)
$$
gives \eqref{eq:radialoperator}, since $V_s=r_{ss}+r_s/(n-k-\varepsilon)$.
All factors divided out are positive, so the calculation is reversible.
The only other second derivatives involving $s$ occur in
$\partial_s\nu$ and $\nabla^\top(1/((\nu\cdot\omega)V))$.
Both depend on $\nabla_{\mathbb S}r_s$ and are independent of $r_{ss}$.
\end{proof}

For the logarithmic radius $v=\log r$, set
\[
\mathscr G_\varepsilon(v)=e^{-v}\mathscr F_\varepsilon(e^v).
\]
Use \eqref{eq:radialoperator} to define $\mathscr G_0$ by setting
$\varepsilon=0$. The round metric and its Laplacian on $\Sphere^{n-1}$
are used for angular derivatives below.

\begin{lemma}[Linearization]\label{star:linearization}
The operator $\mathscr G_\varepsilon$ is invariant under constant
translations of $v$. Its limiting equation $\mathscr G_0(v_0)=0$ holds
along the reference flow. The angular principal matrix of
$D\mathscr G_0(v_0)$ is uniformly positive definite, and the
coefficient of $\partial_s$ is $-1$. For $0<\alpha<1$ and sufficiently
small $\varepsilon>0$, the linearization in a small uniform strip
$C^{2,\alpha}$ neighborhood of $v_0$ is elliptic, with time principal
coefficient $\varepsilon$, mixed coefficients $O(\varepsilon)$,
uniformly positive angular block and uniformly negative time drift.
At any scalar constant $c$, $\mathscr G_\varepsilon(c)=0$ and
\begin{equation}\label{eq:sphericaloperator}
D\mathscr G_\varepsilon(c)
=\varepsilon\partial_{ss}-\partial_s
+\frac{n-k}{(n-1)(n-k-\varepsilon)^2}\Delta_{\mathbb S}
=\mathscr L_{\varepsilon,\infty}.
\end{equation}
\end{lemma}

\begin{proof}
Under dilation of $r$, $V$ has degree one, $h$ and
$1/((\nu\cdot\omega)V)$ have degree minus one, and the tangential
gradient of the latter has degree minus two. Thus every term in
\eqref{eq:radialoperator} has degree one, so
$\mathscr G_\varepsilon(v+t)=\mathscr G_\varepsilon(v)$ for constant
$t$. Its linearization has no zeroth-order term.

Write the physical logarithmic radius as
$f=s/(n-k-\varepsilon)+v$. The normal angle is
$\nu\cdot\omega=(1+|\nabla_{\mathbb S}v|^2)^{-1/2}$, and $V=e^v f_s$.
At $\varepsilon=0$ the equation reads
$$
f_s-k(\nu\cdot\omega)f_s^2
\frac{\sigma_k(e^v h)}{\sigma_{k-1}(e^v h)}=0.
$$
Its positive root is precisely the reference equation. At $v=v_0$,
the coefficient of $\partial_s$ in its derivative is therefore $-1$.
The symmetric shape matrix $e^v h$, in an orthonormal radial frame, is
$$
\frac{(I+\nabla_{\mathbb S}v\otimes\nabla_{\mathbb S}v)^{-1/2}
\bigl(I+\nabla_{\mathbb S}v\otimes\nabla_{\mathbb S}v
-\nabla_{\mathbb S}^2v\bigr)
(I+\nabla_{\mathbb S}v\otimes\nabla_{\mathbb S}v)^{-1/2}}
{\sqrt{1+|\nabla_{\mathbb S}v|^2}}.
$$
Its derivative in an angular Hessian direction is a negative congruence.
Consequently the angular principal matrix of $D\mathscr G_0(v_0)$ is
$$
\begin{aligned}
\frac{k f_s^2}{1+|\nabla_{\mathbb S}v|^2}
&(I+\nabla_{\mathbb S}v\otimes\nabla_{\mathbb S}v)^{-1/2}\\[-2pt]
&\quad\cdot D\left(\frac{\sigma_k}{\sigma_{k-1}}\right)(e^v h)
(I+\nabla_{\mathbb S}v\otimes\nabla_{\mathbb S}v)^{-1/2}
\quad\hbox{at }v=v_0.
\end{aligned}
$$
The derivative of the curvature quotient is positive definite on
$\Gamma_k$. The bounds of Proposition~\ref{star:reference} make the ellipticity
uniform on the entire cylinder.

By Lemma~\ref{star:equation} the time principal coefficient is
$\varepsilon$ and the mixed coefficients are $O(\varepsilon)$.
Smooth dependence on the displayed derivatives preserves the positive
angular block and negative drift in a uniform neighborhood of $v_0$.
The Schur complement there is $\varepsilon-O(\varepsilon^2)>0$.

At a constant $v=c$, $e^v h=I$ and $f_s=1/(n-k-\varepsilon)$.
The residual reduces to $f_s-(n-k-\varepsilon)f_s^2=0$.
Its linearization in the time variables is
$\varepsilon\partial_{ss}-\partial_s$. Rotation invariance and
homogeneity give
$D(\sigma_k/\sigma_{k-1})(I)[M]=(n-k)\operatorname{tr}M/(k(n-1))$.
The normal-angle and tangential-gradient terms have zero first
variation at a constant. Substitution gives \eqref{eq:sphericaloperator}.
\end{proof}

\section{An inverse estimate on the half-cylinder}
\label{smsp:linearinverse}

In this section we prove a weighted inverse estimate for the linearized
equation on the half-cylinder, with inverse norm bounded by
$C\varepsilon^{-6}$ for small $\varepsilon=p-1>0$.
Here $\varepsilon$ is the coefficient of the second time derivative
in Lemma~\ref{star:linearization}. The initial radial graph is fixed,
so the perturbation has zero Dirichlet data at $s=0$. Fix $0<\alpha<1$.
For $m\ge0$ and $\eta>0$ set
$$
\|f\|_{C^{m,\alpha}_\eta}
=\sup_{j=0,1,\ldots}
e^{\eta j}
\|f\|_{C^{m,\alpha}(\mathbb S^{n-1}\times[j,j+1])}.
$$
The norms on each strip use ordinary, isotropic space-time derivatives.
Let $X_\eta$ consist of functions $w$ with $w(\cdot,0)=0$ and
$w-w_\infty\in C^{2,\alpha}_\eta$ for a constant $w_\infty$. Set
$$
\|w\|_{X_\eta}
=|w_\infty|+\|w-w_\infty\|_{C^{2,\alpha}_\eta}.
$$
The two exponents of the constant spherical mode of
$\mathscr L_{\varepsilon,\infty}$ are $0$ and $1/\varepsilon$.
For bounded solutions the forcing and finite boundary value determine
the constant at infinity.

\begin{lemma}[Half-cylinder inverse]\label{smsp:inverse}
Consider a family of operators
$$
\mathscr L_\varepsilon w
=\varepsilon w_{ss}
+2\varepsilon c^{0i}w_{si}
+c^{ij}\nabla_{ij}w
+\beta^0w_s+\beta^i\nabla_iw+\gamma w
$$
on $\mathbb S^{n-1}\times[0,\infty)$. Assume that
$c^{0i},c^{ij},\beta^0,\beta^i,\gamma$ have uniform strip
$C^\alpha$ bounds, the symmetric spatial block $(c^{ij})$ is uniformly
positive definite, $\beta^0\le-1/C<0$, and $\gamma\le0$.
Assume also that the coefficients of the difference
$\mathscr L_\varepsilon-\mathscr L_{\varepsilon,\infty}$ have strip
$C^\alpha$ norms bounded by $Ce^{-\rho j}$ on $j\le s\le j+1$, where
$\mathscr L_{\varepsilon,\infty}$ is given by \eqref{eq:sphericaloperator}. The constants in these
assumptions are independent of $\varepsilon$. If
$$
0<\eta<\min\{\rho,1/(n-k)\},
$$
then, for all sufficiently small positive $\varepsilon$, the equation $\mathscr L_\varepsilon w=f$
for $f\in C^\alpha_\eta$ admits a unique bounded solution with
$w(\omega,0)=0$. This solution belongs to $X_\eta$, and
$$
\|w\|_{X_\eta}
\le C\varepsilon^{-6}\|f\|_{C^\alpha_\eta}.
$$
The constant is uniform in $\varepsilon$ and depends only on the
coefficient bounds, the decay rate and $\eta,\alpha,n,k$.
\end{lemma}

\begin{proof}
First solve the Dirichlet problem with zero data at both ends of
$\mathbb S^{n-1}\times[0,T]$. For small $\varepsilon$ the principal
matrix is positive definite by its Schur complement. The linear elliptic Dirichlet theorem and the maximum
principle apply on this compact smooth cylinder for each fixed
$\varepsilon>0$ (\cite[Chs.~3 and 6]{GilbargTrudinger2001}, in local charts). The nonnegative function $1-e^{-\eta s}$ satisfies
$$
\begin{aligned}
\mathscr L_\varepsilon(1-e^{-\eta s})
&=(-\varepsilon\eta^2+\beta^0\eta)e^{-\eta s}
+\gamma(1-e^{-\eta s})\\
&\le-\frac{\eta}{C}e^{-\eta s}.
\end{aligned}
$$
Multiples of this function bound the solution from both sides by
the maximum principle. The resulting estimate is
$$
\|w\|_{C^0}\le C\|f\|_{C^\alpha_\eta},
$$
independently of $T$ and $\varepsilon$. For fixed $\varepsilon$,
local interior and boundary Schauder compactness as $T\to\infty$
therefore gives a bounded solution on the half-cylinder.

Bounded solutions are unique. Choose a fixed $\delta>0$ and decrease
the upper bound for $\varepsilon$ so that
$\varepsilon\delta<1/C$. Then
$\mathscr L_\varepsilon e^{\delta s}<0$.
For a bounded homogeneous solution with zero data at $s=0$, subtract
an arbitrary positive multiple of $e^{\delta s}$. The difference is
negative at the finite boundary and eventually negative at infinity,
and its strictly positive image under $\mathscr L_\varepsilon$
excludes a positive interior maximum. Letting the multiple tend to
zero, and then applying the same argument to the negative solution,
proves uniqueness.

Next obtain an unweighted second-order estimate. In a sphere chart
near a point of the cylinder use the scales
$$
s=s_0+\varepsilon\tau,
\qquad \theta=\theta_0+\sqrt\varepsilon\,y,
$$
with the second expression interpreted in that chart, and multiply
the equation by $\varepsilon$. The normal principal coefficient
becomes $1$, the mixed coefficients are $O(\sqrt\varepsilon)$, and
the spatial block retains fixed positive bounds. The first- and
zeroth-order coefficients, and all their required $C^\alpha$ norms,
are uniformly bounded. The interior Schauder estimates of \cite[Ch.~6]{GilbargTrudinger2001} apply
on these fixed rescaled charts. At $s=0$ the same argument uses
boundary Schauder estimates with zero Dirichlet data. For the finite
cylinders the identical estimate applies at $s=T$, and its constants
are independent of their length.

Two normal derivatives in the original variables cost
$\varepsilon^{-2}$, and the ordinary H\"older seminorm costs at most
$\varepsilon^{-\alpha}$, giving a total loss of at most
$\varepsilon^{-2-\alpha}$. A finite atlas on the sphere and these
local estimates cover every strip. The rescaled estimates control
nearby H\"older quotients, and derivative suprema control quotients
between more distant points with the same loss.
Since $\alpha<1$, we may use the simpler bound
$$
\|w\|_{C^{2,\alpha}}
\le C\varepsilon^{-3}\|f\|_{C^\alpha_\eta},
$$
where the norm on the left is the uniform unweighted strip norm.
This bound passes to the half-cylinder solution.

Now write its equation as
$$
\mathscr L_{\varepsilon,\infty}w=g,\qquad
g=f+(\mathscr L_{\varepsilon,\infty}
-\mathscr L_\varepsilon)w.
$$
The coefficient decay and the preceding unweighted estimate give
$\|g\|_{C^\alpha_\eta}\le C\varepsilon^{-3}\|f\|_{C^\alpha_\eta}$.

We establish a uniform weighted supremum estimate for the limiting operator. For a spherical harmonic of
degree $\ell$, put
$$
r_\pm(\ell)=\frac{1\pm\sqrt{1+
4\varepsilon\frac{n-k}{(n-1)(n-k-\varepsilon)^2}
\ell(\ell+n-2)}}{2\varepsilon}.
$$
Let $g_\ell$ denote the coefficient of that harmonic in $g$.
The bounded solution with zero data at $s=0$ has coefficient
$$
\begin{aligned}
-\frac1{\varepsilon(r_+-r_-)}\bigg[
&\int_0^s e^{r_-(s-t)}g_\ell(t)\,\mathrm{d}t
+\int_s^\infty e^{r_+(s-t)}g_\ell(t)\,\mathrm{d}t\\
&-e^{r_-s}\int_0^\infty e^{-r_+t}g_\ell(t)\,\mathrm{d}t
\bigg].
\end{aligned}
$$
The derivative of the full-line kernel has jump $1/\varepsilon$
at its pole, which fixes its sign and coefficient. The last term
imposes the finite boundary condition. The positive exponential
solution is excluded by boundedness. For the constant harmonic,
the limiting value is
$$
-\int_0^\infty(1-e^{-t/\varepsilon})g_0(t)\,\mathrm{d}t,
$$
and the difference from this limit is
$$
\int_s^\infty(1-e^{-(t-s)/\varepsilon})g_0(t)\,\mathrm{d}t.
$$
Both the constant and the weighted remainder are uniformly controlled
by $C\|g_0\|_{C^0_\eta}$.

To pass from the modal formulas to a supremum estimate on the sphere, consider
the stable semigroup
$$
\exp\left[
\frac{t}{2\varepsilon}
\left(I-\sqrt{I-
4\varepsilon\frac{n-k}{(n-1)(n-k-\varepsilon)^2}
\Delta_{\mathbb S}}\right)\right].
$$
It is a positive probability mixture of the sphere heat semigroup.
This follows directly from the subordination identity
$$
e^{-t\sqrt{x}}
=\frac{t}{2\sqrt\pi}\int_0^\infty
r^{-3/2}e^{-t^2/(4r)}e^{-rx}\,\mathrm{d}r,
\qquad t,x>0.
$$
The constant harmonic gives total mass one. On zero-mean
functions the heat semigroup has supremum operator norm at most
$C e^{-(n-1)r}$. For $r\le1$ this follows from boundedness.
For $r\ge1$ one uses the spectral gap $n-1$ on $L^2$ for the intermediate time
and bounds the two end factors by fixed-time heat smoothing
from $C^0$ to $L^2$ and from $L^2$ to $C^0$, respectively.
Subordination consequently bounds the stable
semigroup on that subspace by
$$
C e^{-t/(n-k-\varepsilon)}.
$$
The first nonconstant root equals
$r_-(1)=-1/(n-k-\varepsilon)$.
The additional inverse-square-root factor
$1/[\varepsilon(r_+-r_-)]$ in the Green kernel also is a positive
heat mixture, by
$$
x^{-1/2}=\frac1{\sqrt\pi}
\int_0^\infty r^{-1/2}e^{-rx}\,\mathrm{d}r.
$$
It is a contraction and preserves zero mean. The future part of
the Green kernel contains, in addition to the stable semigroup,
$e^{-t/\varepsilon}$. These representations justify summation of
the harmonic formulas as bounded operators on continuous functions.

Since $\eta<1/(n-k)$, convolution of the past kernel with
$e^{-\eta t}$ has bound
$C e^{-\eta s}/(1/(n-k)-\eta)$. The future kernel and the boundary
correction have the same bound uniformly for small $\varepsilon$. Combining the zero-mean and constant parts proves
that the limiting operator maps decaying forcing to a constant plus
a remainder with uniformly bounded $C^0_\eta$ norm. Applying this
estimate to $g$ and using its bound in terms of $f$ gives
$$
|w_\infty|+\|w-w_\infty\|_{C^0_\eta}
\le C\varepsilon^{-3}\|f\|_{C^\alpha_\eta}.
$$

Finally, the equation for $w-w_\infty$ under the original operator
has right-hand side $f-\gamma w_\infty$. The coefficient $\gamma$
decays faster than the chosen weight, so this right-hand side has
the required weighted norm. The finite boundary value is the
constant $-w_\infty$, whose boundary norm is already controlled.
Apply the rescaled Schauder estimate once more on the strips,
now with their weights. The weights vary by bounded factors on
all charts in use. This costs at most another
$\varepsilon^{-3}$ and proves the asserted $X_\eta$ estimate.
\end{proof}

\section{Perturbation of the reference flow}\label{sec:perturbation}

In this section we construct an exact solution of the radial equation
by cancelling the residual along the reference flow through order twelve
and using the inverse estimate of Lemma~\ref{smsp:inverse} in a contraction
argument for the remainder.

\begin{proposition}[Perturbation of the reference flow]\label{star:correction}
Let $v_0=\log r_0$ be the reference solution of
Proposition~\ref{star:reference}, with decay rate $\rho$, and let
$0<\eta<\min\{\rho,1/(n-k)\}$. For every sufficiently small
$\varepsilon>0$ there is a smooth solution on the closed half-cylinder
satisfying
\[
 \mathscr G_\varepsilon(v_\varepsilon)=0,\qquad
 v_\varepsilon(\omega,0)=\log r_\Omega(\omega),\qquad
 \|v_\varepsilon-v_0\|_{X_\eta}\le C\varepsilon .
\]
In particular $v_\varepsilon$ tends to a constant
$v_{\varepsilon,\infty}$ with
$v_\varepsilon-v_{\varepsilon,\infty}\in C^{2,\alpha}_\eta$.
The constant $C$ is independent of $\varepsilon$.
\end{proposition}

\begin{proof}
Choose $\eta<\eta_1<\min\{\rho,1/(n-k)\}$. Seek an approximation
\[
 v_*=v_0+\sum_{j=1}^{12}\varepsilon^jv_j,\qquad
 v_j(\omega,0)=0 .
\]
For $j=1,\ldots,12$, solve successively
\begin{equation}\label{eq:recursion}
 D\mathscr G_0(v_0)v_j=
 -\frac{1}{j!}\left.\frac{\mathrm{d}^j}{\mathrm{d}t^j}
 \mathscr G_t\!\left(v_0+\sum_{i=1}^{j-1}t^i v_i\right)\right|_{t=0},
 \qquad v_j(\cdot,0)=0.
\end{equation}
Adding $t^jv_j$ to the partial sum contributes
$D\mathscr G_0(v_0)v_j$ to the coefficient of $t^j$, so this
choice cancels that order of the residual. By Lemma~\ref{star:linearization} the coefficient of $\partial_s$
in $D\mathscr G_0(v_0)$ is $-1$, the spatial block is uniformly positive
and the zeroth-order coefficient vanishes. Solving for $\partial_s v_j$
therefore gives a forward parabolic equation on the closed sphere.
Smooth initial data and forcing determine a global smooth corrector.

Suppose the forcing and its required derivatives decay with weight
$\eta_1$. Integrating its supremum over time bounds the solution by
the parabolic maximum principle. Local parabolic Schauder estimates (\cite[Ch.~IV]{Lieberman}),
applied on successive strips and at the smooth initial surface,
then bound all required derivatives.
Rewrite the equation as the constant heat operator
$-\partial_s+\Delta_{\mathbb S}/((n-1)(n-k))$ applied to $v_j$,
with the coefficient difference times those bounded derivatives placed
on the right. That difference decays with rate $\rho$. The heat
semigroup on the zero-mean subspace decays with rate $1/(n-k)$.
The tail of the mean is bounded by the integral of the supremum of
this new right-hand side from $s$ to infinity. Thus $v_j$ converges to a constant with error
$O(e^{-\eta_1s})$. Applying the local parabolic estimates to $v_j$
minus this constant, and differentiating the equation, gives the same
decay for the finite number of derivatives required in the recursion.
The estimates also hold on the initial strip. Zero initial values impose no
additional compatibility restriction, since the sphere has no spatial
boundary and the equations determine the initial time derivatives.

To check that \eqref{eq:recursion} produces decaying forcing, observe
that every previous corrector and the reference function tend to constants.
Since $\mathscr G_\varepsilon$ vanishes at every constant for every
small parameter, each constant limiting Taylor coefficient is zero.
Subtracting it leaves at least one decaying factor in every term.
Induction gives the stated decay at each of the twelve orders. Taylor's formula on the fixed
admissible jet neighborhood now gives
$$
\|\mathscr G_\varepsilon(v_*)\|_{C^\alpha_{\eta_1}}
\le C\varepsilon^{13},\qquad
\|v_*-v_0\|_{X_{\eta_1}}\le C\varepsilon.
$$

By Lemma~\ref{star:linearization}, the coefficient of $\partial_{ss}$ in
$D\mathscr G_\varepsilon(v_*)$ is $\varepsilon$,
its mixed coefficients are $O(\varepsilon)$, its spatial ellipticity
is uniform, its time drift is $-1+O(\varepsilon)$, and its
zeroth-order coefficient is zero. Its coefficients converge in strip
$C^\alpha$, with weight $\eta_1$, to the spherical operator \eqref{eq:sphericaloperator}.
The constants are independent of small $\varepsilon$. Lemma~\ref{smsp:inverse} therefore gives an inverse
with zero data at $s=0$ and norm at most $C\varepsilon^{-6}$ from
$C^\alpha_\eta$ to $X_\eta$.

The nonlinear error in
$\mathscr G_\varepsilon(v_*+w)$ after its linear term is removed has
$C^\alpha_\eta$ norm at most $C\|w\|_{X_\eta}^2$, and its difference
at $w,\widetilde w$ is bounded by
$$
C(\|w\|_{X_\eta}+\|\widetilde w\|_{X_\eta})
\|w-\widetilde w\|_{X_\eta}.
$$
This follows from the smooth second-jet dependence and the product
estimates for H\"older norms. The constants in $w$ disappear identically
by translation invariance of $\mathscr G_\varepsilon$, and the nonconstant part decays in the weighted strip norm.
Taylor's formula on a compact neighborhood of the reference jets
then gives the two estimates with constants uniform in $\varepsilon$.

Apply the inverse to the negative of the residual and of this
nonlinear error. On a ball of radius $C\varepsilon^7$ in $X_\eta$,
the residual contributes $O(\varepsilon^7)$ and the nonlinear error
$O(\varepsilon^8)$. The Lipschitz constant is $O(\varepsilon)$.
Fixing the radius constant sufficiently large and then reducing
$\varepsilon$ yields a contraction. If $w$ is its fixed point,
$v_\varepsilon=v_*+w$ has the prescribed data and solves the exact
equation. Moreover,
\[
 \|v_\varepsilon-v_0\|_{X_\eta}
 \le \|v_*-v_0\|_{X_{\eta_1}}+\|w\|_{X_\eta}
 \le C\varepsilon .
\]
For fixed $\varepsilon>0$ the equation is uniformly elliptic near
this solution. Difference quotients and linear Schauder estimates
give interior smoothness. At $s=0$ the smooth Dirichlet data give
the tangential boundary estimates, and the coefficient $\varepsilon$
of $v_{ss}$ allows successive recovery of the normal derivatives.
\end{proof}

\section{Reconstruction and proof of Theorem~\ref*{thm:starshaped}}
\label{sec:main-proof}

In this section we reconstruct the exterior potential from the radial
solution and complete the proof of Theorem~\ref{thm:starshaped} by
verifying \textup{(P1)--(P5)}, applying power comparison, and deriving
the uniform boundary gradient limit.

\begin{proposition}[Reconstruction of the potential]\label{star:reconstruction}
For all sufficiently small $\varepsilon=p-1>0$, let $v_\varepsilon$
be the solution in Proposition~\ref{star:correction} and set
$r_\varepsilon=e^{v_\varepsilon}$. The map
\[
 (\omega,s)\longmapsto
 e^{s/(n-k-\varepsilon)}r_\varepsilon(\omega,s)\omega
\]
is a smooth diffeomorphism of $\Sphere^{n-1}\times[0,\infty)$ onto
$\overline E$. In these coordinates $u_p=e^{-s/\varepsilon}$ solves
\eqref{eq:exterior} and satisfies \textup{(P1)--(P5)} of
Theorem~\ref{thm:starshaped}, with
$\mathsf a_p=e^{\alpha_p v_{\varepsilon,\infty}}$.
\end{proposition}

\begin{proof}
Propositions~\ref{star:reference} and~\ref{star:correction} give
uniform positive lower bounds for $r_\varepsilon$,
$\nu\cdot\omega$, and $V=(r_\varepsilon)_s+
r_\varepsilon/(n-k-\varepsilon)$. Along each ray the physical radius
starts at $r_\Omega$, has derivative $e^{s/(n-k-\varepsilon)}V>0$
and tends to infinity. The inverse function theorem gives the stated
diffeomorphism, including the finite boundary. Lemma~\ref{star:equation}
gives
\[
 Du_p=-\frac{u_p}{\varepsilon}
 \frac{e^{-s/(n-k-\varepsilon)}}{(\nu\cdot\omega)V}\nu\ne0 .
\]
In particular $u_p$ is smooth, equals one on $\partial\Omega$,
lies strictly between zero and one in $E$, and tends to zero at
infinity. Equations \eqref{eq:radialoperator} and
\eqref{eq:newtonidentity} give \eqref{eq:exterior}. The normalized
levels are uniformly $C^2$-close to the reference family, whose
curvature vectors lie in a compact subset of $\Gamma_k$. Hence these
radial graphs, and their physical dilations, are strictly $k$-convex.

To check the branch in dimension $n$, divide $M[u_p]$ by
$(u_p/\varepsilon)e^{-2s/(n-k-\varepsilon)}$. By \eqref{eq:physicalblock}
the tangential block, mixed column, and normal entry of the resulting
matrix are, respectively,
\[
 \frac{h}{(\nu\cdot\omega)V},\qquad
 \sqrt{\frac{\varepsilon}{k}}\,
     \nabla^\top\!\left(\frac1{(\nu\cdot\omega)V}\right),\qquad
 \frac{\varepsilon}{k}e^{2s/(n-k-\varepsilon)}D_\nu|Ds|
       -\frac1{k(\nu\cdot\omega)^2V^2}.
\]
The formula for $D_\nu|Ds|$ in Lemma~\ref{star:equation} and the
space-time $C^2$ bounds give a uniform bound for
$e^{2s/(n-k-\varepsilon)}D_\nu|Ds|$. Along the reference family
the reciprocal of $(\nu_0\cdot\omega)((r_0)_s+r_0/(n-k))$ equals
$k\sigma_k(h_0)/\sigma_{k-1}(h_0)$. Thus, in corresponding
orthonormal frames on the radial graphs, the normalized matrices
differ by $O(\sqrt\varepsilon)$ uniformly from
\[
 k\frac{\sigma_k(h_0)}{\sigma_{k-1}(h_0)}
 \operatorname{diag}\!\left(h_0,-\frac{\sigma_k(h_0)}{\sigma_{k-1}(h_0)}\right).
\]
For $1\le j<k$, the $j$th symmetric polynomial of the block diagonal
factor is
\[
 \sigma_j(h_0)-\frac{\sigma_k(h_0)}{\sigma_{k-1}(h_0)}
                         \sigma_{j-1}(h_0)>0 .
\]
Indeed, Newton's inequalities imply that
$\sigma_j(h_0)/\sigma_{j-1}(h_0)$ strictly decreases with $j$ up
to order $k$. Compactness of the reference curvature range gives a
uniform positive lower bound for each of these expressions. Hence
$\sigma_j(M[u_p])>0$ for $j<k$, while the exact equation gives
$\sigma_k(M[u_p])=0$. Adding $tI$ for small $t>0$ moves the matrix
into $\Gamma_k$, since the derivative of $\sigma_k$ in this direction
is $(n-k+1)\sigma_{k-1}(M[u_p])>0$. This proves
$M[u_p]\in\partial\Gamma_k$, also for $k=1$ with $\sigma_0=1$.
The outward Newton field is smooth by noncriticality and has zero
divergence by \eqref{eq:newtonidentity}, classically and
distributionally. These facts establish \textup{(P1)--(P4)}.

For \textup{(P5)}, invert
\[
 \log|x|=\frac{s}{n-k-\varepsilon}
       +v_\varepsilon(x/|x|,s).
\]
The derivative of the right-hand side with respect to $s$ is
$V/r_\varepsilon$, bounded away from zero.
Proposition~\ref{star:correction} gives
$v_\varepsilon-v_{\varepsilon,\infty}=O(e^{-\eta s})$ with the
same decay for its first two derivatives. The scalar identity first gives
$e^{-\eta s}=O(|x|^{-\eta(n-k-\varepsilon)})$. Differentiating
it once and twice in $(\log|x|,x/|x|)$ shows that the error in
\[
 s=(n-k-\varepsilon)(\log|x|-v_{\varepsilon,\infty})
           +O(|x|^{-\eta(n-k-\varepsilon)})
\]
has the same bound through two derivatives in these coordinates.
Since each Euclidean derivative contributes a factor $|x|^{-1}$,
$u_p=e^{-s/\varepsilon}$ satisfies, uniformly in direction,
\[
 D^j\!\left(u_p(x)-e^{\alpha_p v_{\varepsilon,\infty}}|x|^{-\alpha_p}\right)
 =O_p\!\left(|x|^{-\alpha_p-j-\eta(n-k-\varepsilon)}\right),
 \qquad j=0,1,2 .
\]
This estimate is for each fixed $p$.
\end{proof}

\begin{proof}[Proof of Theorem~\ref{thm:starshaped}]
Choose $0<\varepsilon_0<n-k-1$ so that Propositions~\ref{star:correction}
and~\ref{star:reconstruction} apply whenever $0<\varepsilon=p-1<\varepsilon_0$.
The reconstructed potential $u_p=e^{-s/\varepsilon}$ satisfies
\textup{(P1)--(P5)}. In particular, $Du_p$ never vanishes,
$M[u_p]\in\partial\Gamma_k$, and the levels are strictly $k$-convex.
Proposition~\ref{star:boundary} therefore gives uniqueness on the stated
viscosity branch and the interior-ball comparison.

At $s=0$ the surface is $\partial\Omega$ and
$r_\varepsilon(\omega,0)=r_\Omega(\omega)$. Let $\nu$ be its outward
unit normal at $r_\Omega(\omega)\omega$, and let $x(\omega,s)$ denote the
radial parametrization of Proposition~\ref{star:reconstruction}.
Since $Ds$ is normal to the level and $Ds\cdot x_s=1$, we have
$|Ds|=(\nu\cdot x_s)^{-1}$ there.
Differentiating $u_p=e^{-s/\varepsilon}$ gives
\[
 \varepsilon|Du_p|
 =\frac{1}{(\nu\cdot\omega)
         \bigl((r_\varepsilon)_s+r_\Omega/(n-k-\varepsilon)\bigr)} .
\]
Let $r_0$ be the normalized radius of the reference flow from the same
boundary (Proposition~\ref{star:reference}). By Proposition~\ref{star:correction},
$(r_\varepsilon)_s(\cdot,0)\to(r_0)_s(\cdot,0)$ uniformly.
The limiting denominator is the outward normal speed of the reference
flow at $s=0$, namely $\sigma_{k-1}(h_\delta)/(k\sigma_k(h_\delta))$.
It has a positive minimum on the compact boundary, so taking reciprocals
gives the uniform limit
$k\sigma_k(h_\delta)/\sigma_{k-1}(h_\delta)=(n-k)H_k/H_{k-1}$,
which proves \eqref{eq:boundarylimit}.
\end{proof}

\bibliographystyle{plain}
\bibliography{references/paper2/reference}
\end{document}